\documentclass[12pt,a4paper,reqno]{amsart}
\usepackage[all]{xy}
\usepackage{amsmath}
\usepackage{amsfonts}
\usepackage{amssymb}
\usepackage{mathrsfs}
\usepackage{calligra}
\usepackage{tensor}
\usepackage[utf8]{inputenc}
\usepackage{xcolor}
\usepackage{graphicx}

\theoremstyle{plain}
\newtheorem{theorem}{Theorem}[subsection]
\newtheorem{mytheorem}{Theorem}
\newtheorem{remark}[theorem]{Remark}
\newtheorem{lemma}[theorem]{Lemma}
\newtheorem{proposition}[theorem]{Proposition}
\newtheorem{corollary}[theorem]{Corollary}
\newtheorem{definition}[theorem]{Definition}

\newcommand{\Bmu}{\mbox{$\raisebox{-0.59ex}
  {$l$}\hspace{-0.18em}\mu\hspace{-0.88em}\raisebox{-0.98ex}{\scalebox{2}
  {$\color{white}.$}}\hspace{-0.416em}\raisebox{+0.88ex}
  {$\color{white}.$}\hspace{0.46em}$}{}}

\DeclareMathAlphabet{\mathbbold}{U}{bbold}{m}{n}

\DeclareFontEncoding{OT2}{}{} 

\newcommand{\scHom}{\mathscr{H}\text{\kern -3pt {\calligra\large om}}}
\newcommand{\scExt}{\mathscr{E}\text{\kern -3pt {\calligra\large xt}}}

\def\O{{\mathcal O}}
\def\P{{\mathcal P}}
\def\Q{{\mathbb Q}}

\def\V{{\mathcal V}}

\def\fY{{\mathfrak Y}}
\def\Z{{\mathbb Z}}

\DeclareMathOperator{\Div}{Div}
\DeclareMathOperator{\coh}{H}
\DeclareMathOperator{\Gal}{Gal}

\DeclareMathOperator{\image}{Im}
\DeclareMathOperator{\Nm}{N}

\DeclareMathOperator{\Hom}{Hom}
\DeclareMathOperator{\ord}{ord}

\begin{document}
\title[capitulation kernels]
{capitulation kernels of class groups over $\Z_p^d$-extensions}

\author{King Fai Lai}
\address{School of Mathematics and Statistics\\
Henan University \\
 Jin Ming Avenue, Kaifeng, Henan, 475004, China}
\email{kinglaihonkon@163.com}

\author[Tan]{Ki-Seng Tan$^\dag$}
\address{Department of Mathematics\\
National Taiwan University\\
Taipei 10764, Taiwan} \email{tan@math.ntu.edu.tw}

\thanks{$^{\dag}$ The author was supported in part by the Ministry of
Science and Technology of Taiwan, MOST 103-2115-M-002 -008 -MY2,
MOST 105-2115-M-002 -009 -MY2.\\
It is our pleasure to thank NCTS/TPE for supporting a number of meetings of the authors in National Taiwan University.
}



\maketitle

\begin{abstract}
We generalize Iwasawa's theorems on class groups over $\Z_p$-extensions to all $\Z_p^d$-extensions.
\end{abstract}

\section{Introduction}\label{s:aim} Let $k$ be a number field and let $K/k$ be an abelian extension.
If $K/k$ is a  finite extension,
we define a map from the group of ideals of $k$ to that of $K$
 by taking $\mathfrak a$ to the fractional ideal $\mathfrak a \O_{K}$.
This map induces the homomorphism 
 $c_{K/k}: C_k\longrightarrow  C_K$ of class groups which is called the
capitulation map. Previous studies indicate that there is certain connection between the order of $\ker (c_{K/k})$ and the ramification of $K/k$.
A famous theorem of Suzuki  \cite{suz91} that generalizes both Hilbert's theorem 94 and
the principal ideal theorem of class field theory
says if $K/k$ is unramified, then the
order of $\ker(c_{K/k})$ is 
divisible by the degree $[K:k]$. 
More results in this aspect can be found in \cite{gon07}.
In contrast,
Iwasawa \cite{iwa73} proves that if $K/k$ is a $\Z_p$-extension, so it is almost totally ramified at some place dividing $p$,
then the 
the capitulation kernel $\dot X_K$ (see below) 
is pseudo-null.
In this note, we generalize this result of Iwasawa to every $\Z_p^d$-extension of $k$, and then use it to establish a pseudo-isomorphism
of Iwasawa modules that generalizes \cite[Theorem 11]{iwa73}.

From now on let $K/k$ be a $\Z_p^d$-extension. 
For a number field $E\subset K$, let $A_E$ denote the Sylow $p$-subgroup
of the class group ${C}_E$ of $E$ and let $A'_E$ be the quotient of $A_E$ modulo the subgroup generated by ideals above $p$.
Define 
$$\dot{A}_E=\bigcup_{E\subset E'\subset K}\ker(c_{E'/E}),\;\;\dot{A}'_E=\bigcup_{E\subset E'\subset K}\ker(c'_{E'/E}),$$
where $E'/E$ runs through finite sub-extensions in $K/E$ and 
$$c'_{E'/E}:A'_E\longrightarrow A'_{E'}$$
is induced from $c_{E/E'}$. Note that $\dot A_{E}$ is a subgroup of $A_E$. Put 
$$\dot{X}_K=\varprojlim_{k\subset E\subset K}\dot{A}_E,\;\;\dot{X}'_K=\varprojlim_{k\subset E\subset K}\dot{A}'_E,$$
with the limit taken over the norm maps.
Let $\sim$ denote pseudo-isomorphism.
\begin{mytheorem}\label{t:x0}
We have $\dot{X}_K\sim 0$ and $\dot{X}'_K\sim 0$.
\end{mytheorem}
Theorem \ref{t:x0}
 is proved in \S\ref{su:x0} by following the path of Iwasawa \cite[Theorem 10]{iwa73}. 
 For $d\geq 2$, the Iwasawa algebra is much more
 complicated than that of a $\Z_p$-extension, this might be the reason why the content of the theorem has seldom discussed since
 the publication of Iwasawa's paper.  The breakthrough in our proof is the duality in Lemma \ref{l:monsky}(c). Having it, we use
 Monsky's theorem \cite{mon81} (see \S\ref{su:monsky}) to show Proposition \ref{p:monsky}, which then leads to the proof of the theorem.


Write $\Gamma$  for the Galois group of $K/k$ and $\Lambda_\Gamma$ for the
Iwasawa algebra $\Z_p[[\Gamma]]$. 
Define the $\Lambda_\Gamma$-modules
$$W_K:=\varprojlim_{k\subset E\subset K} \Hom(A_{E},\Q_p/\Z_p),\;\;W_K':=\varprojlim_{k\subset E\subset K}\Hom(A'_{E},\Q_p/\Z_p),$$ 
with the limits taken over the homomorphisms dual to the capitulation maps,
and 
$$X_K:=\varprojlim_{k\subset E\subset K} A_E,\;\;X'_K:=\varprojlim_{k\subset E\subset K} A'_E.$$

The following theorem might be known to experts, for the convenience of the readers, we include
a proof of it in \S\ref{su:ml}.
\begin{mytheorem}\label{t:m}
The $\Lambda_\Gamma$-modules $W_K$ and $W'_K$ are finitely generated.
\end{mytheorem}

Let
$${}^\sharp:\Lambda_\Gamma\longrightarrow \Lambda_\Gamma$$
be the involution of $\Z_p$-algebra induced from $\Gamma\longrightarrow\Gamma$, $\gamma\mapsto\gamma^{-1}$.
For a $\Lambda_\Gamma$-module $\mathfrak D$, let
$\mathfrak D^\sharp$ denote the module with $\mathfrak D$ as the underlying $\Z_p$-module while $\Gamma$ acts on $\mathfrak D$
via ${}^\sharp:\Gamma\longrightarrow \Gamma$.

\begin{mytheorem}\label{t:mw}  
We have
$$W_K\sim X_K^\sharp,\;\;W_K'\sim X_K'^\sharp.$$
\end{mytheorem}
Suppose $k_0$ is a subfield of $k$ such that $K/k_0$ is an abelian extension having Galois group
$\Gal(K/k_0)=\Gamma\times \Theta$ with $\Theta$ finite of order prime to $p$. Then every pro-$p$ $\Gal(K/k_0)$-module
$\mathfrak D$ can be written as
$$\mathfrak D=\bigoplus_{\chi} \mathfrak D_\chi,$$
where $\chi$ runs through all ${\bar{\Q}_p}^*$-valued characters of $\Theta$ and
$$\mathfrak D_\chi=\{x\in\mathfrak D\;\mid\; \tensor[^\sigma] x{}=\chi(\sigma)\cdot x,\;\text{for all}\;\sigma\in\Theta\}.$$
It is clear that $(\mathfrak D^\sharp)_\chi=(\mathfrak D_{\chi^{-1}})^\sharp$, and we have (see the last paragraph of \S\ref{su:ml})
\begin{equation}\label{e:wchi}
W_\chi\sim (X_{\chi^{-1}})^\sharp,\;\;W'_\chi\sim X_{\chi^{-1}}'^\sharp.
\end{equation}
These pseudo-isomorphisms imply the equalities of characteristic ideals
\begin{equation}\label{e:characteristicideal}
\mathrm{CH}_\Gamma (W_\chi)=\mathrm{CH}_\Gamma (X_{\chi^{-1}}^\sharp),\;\;\mathrm{CH}_\Gamma (W_\chi')=\mathrm{CH}_\Gamma (X_{\chi^{-1}}'^\sharp).
\end{equation}

Theorem \ref{t:mw}, proved in \S\ref{su:ml}, is not entirely new, when $d=1$, it is by Iwasawa
\cite[Theorem 11]{iwa73}. Suppose $K/k$ has ramification locus $S$. 
Nekov\'{a}\v r defines in \cite[\S 9.5]{nek06} a morphism from $X_K$ to the Iwasawa adjoint $E^1(W_K^\sharp)$
of $W_K^\sharp$ (they are pseudo-isomorphic) and shows the cokernel is pseudo-null under the condition:
\vskip5pt
\noindent
(Dec 2): For every $v\in S$, the decomposition group $\Gamma_v\simeq \Z_p^{r(v)}$, with $r(v)\geq 2$.
\vskip5pt 
A similar result for $X_K'$ and $W_K'^\sharp$ is given in \cite[\S 9.4]{nek06}, while Vauclair, in collaboration with Nekov\'{a}\v r, 
proves in \cite[Theorem 7.6]{vau09}
that there is a natural morphism $\alpha:W_K'\longrightarrow E^1(X_K'^\sharp)$ with pseudo-isomorphic kernel and cokernel, and if (Dec 2) holds then $\alpha$ is a pseudo-isomorphism. This also implies unconditionally the equality of characteristic ideals as in \eqref{e:characteristicideal} (see \cite[Corollary 6.8]{vau09})
$$\mathrm{CH}_\Gamma (W'_K)=\mathrm{CH}_\Gamma (X_K'^\sharp).$$

Our proof of Theorem \ref{t:mw} uses the theory of $\Gamma$-system established in \cite{lltt18}. It turns out Theorem \ref{t:x0} is
the key ingredient for having the theory work in our situation. Indeed, a sensible homomorphism $X_K\longrightarrow W_K^\sharp$
should be compatible with the duality between $A_E$ and $\Hom(A_E,\Q_p/\Z_p)$ at all finite layers, and hence must factors through 
the quotient of $X_K$ modulo $\dot X_K$. Therefore, the pseudo isomorphism between $X_K$ and $W_K^\sharp$ cannot established 
without having $\dot X_K\sim 0$, the same with $X'_K$ and $W_K'^\sharp$.

It is worthwhile to mention that over global function fields of characteristic $p$
our theorems hold trivially, because the global unit group $U_K$ has trivial $p$-primary part (cf. \S \ref{su:exact} and \S \ref{su:ml}, especially the exact sequence
\eqref{e:coh1}, Lemma \ref{l:equiv} and the proof of Theorem \ref{t:mw}).

We thank I. Longhi for helping us with the proof of Lemma \ref{l:inverselim},
we also thank D. Vauclair for a communication on our and his work.

\medskip

\subsection{Notation and preliminary remarks}\label{su:notation}

Let $E$ stand for a finite intermediate extension of $K/k$ and denote $\Gamma_E=\Gal(K/E)$.
Put $\Gamma^{(n)}:=\Gamma^{p^n}$, $k_n:=K^{\Gamma^{(n)}}$, the $n$th layer of $K/k$.
Let $\mathscr I_\Gamma$ denote
the augmentation ideal of $\Lambda_\Gamma$ and put $\mathscr I_n:=\ker (\Lambda_\Gamma \longrightarrow \Z_p[\Gal(k_n/k)])$ so that $\mathscr I_0=\mathscr I_\Gamma$.


For a set $T$ of places of $k$, let
$T_E$ denote the set of places of $E$ sitting over $T$. Let $\P$ be the set of places of $k$ above $p$.
For a Galois group $G$ of number field extension, let $G_w$ and $G_w^0$ denote respectively the decomposition subgroup and
the inertia subgroup at a place $w$.

For a topological group $B$, let $B^\vee$ denote  $\Hom(B,\Q_p/\Z_p)$, where $\Q_p/\Z_p$ is endowed with the discrete topology.
In most cases, $B$
is pro-$p$ or $p$-primary and discrete so that $B^\vee$ is its Pontryagin dual. We also identify $B^\vee$ with $\coh^1(B,\Q_p/\Z_p)$
when $B$ is a compact group.
For a functor $\lozenge$ on the category of finite intermediate extension of $K/k$, we abbreviate $\lozenge_{k_n}=\lozenge_n$.
If $\mathfrak Q \longrightarrow \mathfrak D$ is a morphism of two such functotors
with $\mathfrak Q_E$, $\mathfrak D_E$ compact and 
$\mathfrak Q_E\longrightarrow \mathfrak D_E$
surjective for every $E$, then $\mathfrak Q_K\longrightarrow \mathfrak D_K$ is also surjective. To see this, we apply the dual
$\xymatrix{\mathfrak D_E^\vee \ar@{^{(}->}[r] & \mathfrak Q_E^\vee}$, which implies the injectivity of
$\xymatrix{\mathfrak D_K^\vee \ar@{^{(}->}[r] & \mathfrak Q_K^\vee}$, and then obtain the desired surjectivity via the duality \cite{kap50}.

Since $X_K$ is torsion, it is annihilated by some non-zero $f\in\Lambda_\Gamma$.
Then $W_K$ is annihilated by $f^\sharp$, whence torsion.
To see this, we may assume that for every intermediate extension $E$ of $K/k$, the extension
$K/E$ contains no nontrivial unramified intermediate extension.
Then the restriction map $\Hom(A_E,\Q_p/\Z_p) \longrightarrow \Hom(A_{E'},\Q_p/\Z_p)$ is injective for $E\subset E'\subset K$, and hence
$X_K\longrightarrow A_E$ is surjective.
This implies $f\cdot A_E=0$, so
$f^\sharp\cdot \Hom(A_E,\Q_p/\Z_p)=0$ and $f^\sharp\cdot W_K=0$.


\section{Towers of $\Gamma$-modules}\label{s:monsky}

\subsection{Monsky's Theorem}\label{su:monsky}

Endow $\Bmu_{p^\infty}$ with the discrete topology. Let $\hat\Gamma$ denote all continuous characters from $\Gamma$ to $\Bmu_{p^\infty}$.
Every $\chi\in\hat\Gamma$ is of finite order, it factors through $\Gamma\longrightarrow \Gamma/\Gamma^{(n)}$, for some $n$,
and hence extends uniquely to a continuous $\Z_p$-algebra homomorphism $\chi:\Lambda_\Gamma\longrightarrow \O$, where $\O$
is the ring of integers of $\bar\Q_p$.

\begin{definition}
\begin{enumerate}
\item[(1)] Define the {\em{zero set}} of an element $\theta\in\Lambda_\Gamma$ to be
$$\Delta_\theta:=\{\chi\in\hat\Gamma\;\mid\;\chi(\theta)=0\}.$$
\item[(2)]
A $\Z_p$-flat $Z$ of codimension $m$ is a subset of $\hat\Gamma$ consisting of solutions $\chi$ to
the system of equations
\begin{equation}\label{e:def}
\chi(\xi_{j})=\zeta_{j},\;j=1,...,m,
\end{equation}
where $\xi_{1},...,\xi_{m}$ are elements of $\Gamma$, extendable to a $\Z_p$-basis, 
$\zeta_{1},...,\zeta_{m}\in\Bmu_{p^\infty}$. 
\\
\item[(3)]  We say a set $\tau_1,....\tau_c$ of $\Z_p$-generators of $\Gamma$ is {\em{tight}}, if each $\tau_i\notin \Gamma^{(1)}$ and the topological closure of $<\tau_i>$ and $<\tau_j>$ are  distinct for $i\not=j$.
\end{enumerate}
\end{definition}

\begin{theorem}[{Monsky, \cite[Theorem 2.2,2.6]{mon81}} ]\label{t:monsky} The zero set $\Delta_f$ of a non-zero $f\in\Lambda_\Gamma$ is a proper subset of $\hat\Gamma$. There are $\Z_p$-flats $Z_1,...,Z_l$ such that
$$\Delta_f=Z_1\cup\cdots\cup Z_l.$$
\end{theorem}

Let us introduce more notation. For $n\geq 0$, $\sigma\in\Gamma$, $\sigma\not=id$,
denote
$$\omega_{\sigma,n}:=\sigma^{p^n}-1,$$
which is regarded as an element of $\Lambda_\Gamma$. Put $\omega_{\sigma,-1}=1$. For $m\geq n\geq -1$, write
$$\nu_{\sigma,n,m}:=\omega_{\sigma,m}/\omega_{\sigma,n}$$
which is also an element of $\Lambda_\Gamma$.
Because $\nu_{\sigma,r,r+1}$ is the $p^{r+1}$th cyclotomic polynomial in $\sigma$, it is irreducible in $\Z_p[\sigma]$.
Hence in $\Lambda_\Gamma$, if $r\not=r'$, then $\nu_{\sigma,r',r'+1}$ is relatively prime to
$\nu_{\sigma,r,r+1}$. Also, if $\sigma$ and $\sigma'$ are linearly independent over $\Z_p$, then $\nu_{\sigma,n,m}$
and $\nu_{\sigma',n',m'}$ are relatively prime in $\Lambda_\Gamma$.

Due to technical reason, we will need to deal with ideals other than $\mathscr I_n$, 
especially when $d\geq 2$. We introduce the ideals 
$$\mathscr{J}_n:=(\nu_{\tau_1,0,n},...,\nu_{\tau_c,0,n})\subset \Lambda_\Gamma,$$
where $\tau_1,....\tau_c$ is a chosen set of {\em{tight}}  $\Z_p$-generators of $\Gamma$. It is clear that for $m\geq n$, the inclusion $\mathscr{J}_m\subset \mathscr{J}_n$ holds.
Since $\omega_{\tau,0}\cdot\nu_{\tau,0,n}=\omega_{\tau,n}$, the ideal
$\mathscr{I}_n=(\omega_{\tau_1,n},...,\omega_{\tau_c,n})$ is inside $\mathscr{J}_n$.
For a chosen $\Z_p$-basis $\sigma_1,...,\sigma_d$ of $\Gamma$, we also need to consider ideals,
$$\mathscr{I}_{\underline r,\underline n}:=(\nu_{\sigma_1,r_1,n_1},\nu_{\sigma_2,r_2,n_2},...,\nu_{\sigma_d,r_d,n_d})
\subset \Lambda_\Gamma,$$
with $\underline n=(n_1,...,n_d)$, $\underline r=(r_1,...,r_d)\in\Z^d$ such that $r_i\geq -1$ and $n_i>r_i$, for every $i$,
the latter condition will be abbreviated as $\underline n>\underline r$.

For an ideal $I=(\theta_1,...,\theta_l)\subset \Lambda_\Gamma$, denote
$$\Delta_{I}:=\bigcap_{\theta\in I} \Delta_\theta=\bigcap_{i=1}^l \Delta_{\theta_i}.$$
Write $\Delta_{\underline r,\underline n}$ for
$\Delta_{\mathscr{I}_{\underline r,\underline n}}$.
If in \eqref{e:def} each $\zeta_j$ is of order $p^{r_{j}}$, then
\begin{equation}\label{e:Z}
Z\subset \Delta_{(\nu_{\xi_1,r_1-1,r_1}, \cdots,\nu_{\xi_m,r_m-1,r_m})}=
\Delta_{\nu_{\xi_1,r_1-1,r_1}}\cap \cdots \cap \Delta_{\nu_{\xi_m,r_m-1,r_m}}.
\end{equation}

\begin{lemma}\label{l:monsky} The following statements hold true:
\begin{enumerate}
\item[(a)] If the order of $\chi\in\hat\Gamma$ does not exceed $p^n$ or $p^{n_i}$, for every $i=1,...,d$, then for all
$x\in\mathscr{I}_{\underline r,\underline n}+\mathscr{J}_n$ the value $\chi(x)$ is divisible by the minimum of
$p^n$ and $p^{n_1-r_1},...,p^{n_d-r_d}$.
\item[(b)] For a fixed $\underline r$, the intersection $\displaystyle{\bigcap_{\underline n>\underline r, n>0}
(\mathscr{I}_{\underline r,\underline n}+\mathscr{J}_n})=0$.
\item[(c)] Let $x\in\Lambda_\Gamma$. Then $x\in\mathscr{I}_{\underline r,\underline n}$ if and only if $\chi(x)=0$ for all
$\chi\in\Delta_{\underline r,\underline n}$.
\end{enumerate}
\end{lemma}
\begin{proof} Let $p^\alpha$ denote the order of $\chi(\sigma)$. For $m\geq \alpha$,
we have $\chi(\nu_{\sigma,r,m})=0$, if $\alpha >r$, while $\chi(\nu_{\sigma,r,m})=p^{m-r}$, if $\alpha\leq r$.
The assertion (a) follows. 

To prove (b), let $x$ be an element of the intersection in question. In view of Theorem \ref{t:monsky}, we have to show
$\chi(x)=0$ for all $\chi\in\Hat\Gamma$.
Let $\alpha$ denote the maximum of $r_1,...,r_d$ . For an integer $\beta>\alpha$, put $\underline n=(\beta,...,\beta)$, $n=\beta$.
Then by (a), the value $\chi(x)$ is divisible by $p^{\beta-\alpha}$. Since this holds for all $\beta$, we must have $\chi(x)=0$.

We prove (c) by induction on $d$. The statement trivially holds for $d=0$. Suppose $d>0$, let $\Gamma'\subset \Gamma$ be the
subgroup topologically generated by elements of the basis $\sigma_1,...,\sigma_d$ other than $\sigma_1$.
If $d=1$, then both $\Gamma'$ and $\hat{\Gamma'}$ are the trivial group and $\Lambda_{\Gamma'}=\Z_p$. In this case,
put $\mathscr{I}_{\underline{r'},\underline{n'}}=(0)$ and $\Delta_{\underline{r'},\underline{n'}}=\hat{\Gamma'}$.
If $d\geq 2$, let $\underline{r'}=(r_2,...,r_d)$ and $\underline{n'}=(n_2,...,n_d)$. Put $T:=\sigma_1-1\in\Lambda_\Gamma$.
Then $\nu_{\sigma_1,r_1,n_1}$ is a distinguished polynomial in $T$ of degree $\delta:=p^{n_1}-[p^{r_1}]$, where $[\cdot]$ is the
Gauss symbol.
The natural map
$\hat\Gamma\longrightarrow\hat{\Gamma'}$, $\chi\mapsto\bar\chi$, induces a surjection
$\rho:\Delta_{\underline r,\underline n}\longrightarrow \Delta_{\underline{r'},\underline{n'}}$ such that the fibre $\rho^{-1}(\bar\chi)$
for each $\bar\chi\in \Delta_{\underline{r'},\underline{n'}}$ is of cardinality $\delta$ and a character $\chi\in\rho^{-1}(\bar\chi)$ is determined
by the value $\chi(\sigma_1)$, or equivalently, by $\chi(T)$. By the Weierstrass division theorem \cite[VII, \S 3, Proposition 5]{bou}, we may assume that
$$x=\sum_{i=0}^{\delta-1} y_i\cdot T^i,  \;\; y_0,...,y_{\delta-1}\in\Lambda_{\Gamma'}.$$
Thus, for $\chi\in \rho^{-1}(\bar\chi)$,
$$0=\chi(x)=\sum_{i=0}^{\delta-1} \bar\chi(y_i)\cdot \chi(T)^i.$$
This means the polynomial $\sum_{i=0}^{\delta-1} \bar\chi(y_i)\cdot X^i$ in $X$ has $\delta$ distinct roots, hence must be trivial.
Therefore, each $y_i$ is annihilated by characters in $\Delta_{\underline{r'},\underline{n'}}$.
The induction hypothesis implies $y_i\in \mathscr{I}_{\underline{r'},\underline{n'}}$ for every $i$, so $x\in\mathscr{I}_{\underline{r},\underline{n}}$.
\end{proof}

\begin{lemma}\label{l:inverselim}
Suppose $\fY$ is a finitely generated $\Lambda_\Gamma$-module. Then
$\bigcap_n\mathscr{J}_n\fY=0$ so that the natural maps $\fY\longrightarrow\fY/\mathscr{J}_n\fY$ induce the isomorphism
$$\xymatrix{\fY\ar[r]^-\sim & \bar\fY:=\varprojlim_{n} \fY/\mathscr{J}_n\fY}.$$
\end{lemma}

\begin{proof} The homomorphism $\fY\longrightarrow\bar\fY$ has
dense image, and  it is surjective because both $\fY$ and $\bar\fY$ are compact. It is sufficient to show
$\bigcap_n\mathscr{J}_n\fY=0$. Lemma \ref{l:monsky} (b) says the assertion holds for $\fY=\Lambda_\Gamma$.
In general, we have a surjective homomorphism $\phi:\bigoplus_{i=1}^l \Lambda_\Gamma\longrightarrow \fY$,
with $\mathscr{J}_n\fY=\phi(\bigoplus_{i=1}^l \mathscr{J}_n)$. Suppose $x\in \bigcap_n\mathscr{J}_n\fY$.
Denote $X=\phi^{-1}(x)$ which is a compact subset of $\bigoplus_{i=1}^l \Lambda_\Gamma$.
Since $X\cap \bigoplus_{i=1}^l \mathscr{J}_n\not=\emptyset$, for every $n$, and $\bigcap_n \bigoplus_{i=1}^l \mathscr{J}_n=0$,
we must have $0\in X$. Therefore, $x=\phi(0)=0$.
\end{proof}

\subsection{The norm maps}\label{su:norm}
Fix a basis $\sigma_1,...,\sigma_d$ of $\Gamma$ and set
$$\nu_{n,m}:=\prod_{i=1}^d \nu_{\sigma_i,n,m}\in \Lambda_\Gamma.$$
Since $\mathscr{I}_n$ is gernerated by $\omega_{\sigma_i,n}$ and $\nu_{\sigma_i,n,m}\cdot \omega_{\sigma_i,n}=\omega_{\sigma_i,m}$, we have
$\nu_{n,m}\cdot\mathscr{I}_n\subset \mathscr{I}_m$.
For a $\Lambda_\Gamma$-module $\fY$, the map sending $y\in \fY$ to $\nu_{n,m}\cdot y$
induces an endomorphism $\xymatrix{\fY/\mathscr{I}_m\fY\ar[r]^-{\nu_{n,m}} & \fY/\mathscr{I}_m\fY}$,
which turns out to be the norm  map of $\Gamma^{(n)}/\Gamma^{(m)}$ acting on $\fY/\mathscr{I}_m\fY$, hence independent of the choice of the basis. 

In general, if for each $n$, there associates a
$\Lambda_\Gamma$ submodule
$\mathfrak I_n\subset \fY$ containing $\mathscr{I}_n\fY$ such that $\nu_{n,m}\cdot \mathfrak I_n\subset \mathfrak I_m$ for $m\geq n$, we shall
also denote the $\Lambda_\Gamma$-map induced from $\nu_{n,m}$ as
$$\xymatrix{\fY/\mathfrak I_n \ar[r]^-{\nu_{n,m}} & \fY/\mathfrak I_m}.$$

Let $\mathscr{J}_n$ be defined by choosing a tight set $\tau_1,...,\tau_c$ of generators of $\Gamma$.
\begin{lemma}\label{l:kernorm} We have $\nu_{n,m}\cdot \mathscr{J}_n\subset \mathscr{J}_m$.
\end{lemma}
\begin{proof} 
We have to show that $\nu_{n,m}\cdot\nu_{\tau_i,0,n}\in \mathscr{J}_m$. Since modulo $ \mathscr{I}_m$, the image
$\nu_{n,m}\cdot \nu_{\tau_i,0,n}$ is
independent of the choice of the basis, we may assume that  $\tau_i=\sigma_1$. Then
$\nu_{n,m}(\nu_{\tau_i,0,n})=\nu_{\sigma_1,0,m}\cdot \nu_{\sigma_2,n,m}\cdot\cdots\cdot \nu_{\sigma_d,n,m}\in\mathscr{J}_m$.
\end{proof}

\begin{definition}\label{d:ddot}For a $\Lambda_\Gamma$-module $\fY$, define
$$\Ddot\fY_{(n,m)}:=\ker(\,\xymatrix{\fY/\mathscr{J}_n\fY \ar[r]^-{\nu_{n,m}} & \fY/\mathscr{J}_m\fY}\,),$$
$$\Ddot\fY_{(n)}:=
\bigcup_{m\geq n} \Ddot\fY_{(n,m)},$$
and
$$\Ddot\fY:=\varprojlim_{n}\Ddot\fY_{(n)}.$$
\end{definition}

\begin{proposition}\label{p:monsky} Suppose $\fY$ is a finitely generated $\Lambda_\Gamma$-module. Then $\Ddot\fY\sim0$.
\end{proposition}
\begin{proof} Consider an exact sequence
$$\xymatrix{0 \ar[r]  & \mathfrak N' \ar[r]  & \fY \ar[r] & \mathfrak Z  \ar[r] & \mathfrak N \ar[r] & 0},$$
with $ \mathfrak N$, $ \mathfrak N'$ pseudo-null and $\mathfrak Z$ a direct sum $\oplus_{i=1}^r\Lambda_\Gamma/(f_i)$ where $f_i$ is either $0$ or a power of an irreducible element of $\Lambda_\Gamma$. Since $ \mathfrak N\oplus \mathfrak N'$ is pseudo-null,
there are relatively prime $g_1,g_2\in\Lambda_\Gamma$ annihilating both $ \mathfrak N$ and $ \mathfrak N'$. These would lead to an exact seuqence
$$\xymatrix{0 \ar[r]  & \mathfrak{a} \ar[r]  & \Ddot\fY \ar[r] & \Ddot{\mathfrak{Z}}  \ar[r] & \mathfrak{b} \ar[r] & 0},$$
with both $\mathfrak a$ and $\mathfrak b$ annihilated by $g_1^2$ and $g_2^2$, hence pseudo-null. Thus, it remains 
to show that $\Ddot\fY=0$ for 
 $\fY=\Lambda_\Gamma/(f)$. 

Let $\bar x\in\Ddot\fY \subset \fY$ (Lemma \ref{l:inverselim}) and let $x\in\Lambda_\Gamma$ be a lifting of $\bar x$.
For each $n$, there is $m>n$ such that
\begin{equation}\label{e:x}
\nu_{n,m}\cdot x\equiv a\cdot f \pmod{\mathscr{J}_m},
\end{equation}
for some $a\in\Lambda_\Gamma$. Put $\omega:=\omega_{\tau_1,0}\cdot\cdots\cdot\omega_{\tau_c,0}$, $y:=\omega x$. Since
$\omega\cdot\mathscr{J}_m\subset \mathscr{I}_m$,
\begin{equation}\label{e:omegay}
\nu_{n,m}\cdot y\equiv \omega a\cdot f \pmod{\mathscr{I}_m}.
\end{equation}

Suppose $f=0$. For every $\chi\in\Delta_{\mathscr{I}_n}$, we have $\chi(\nu_{n,m})\not=0$, while $\chi(\mathscr{I}_m)=0$, so $\chi(y)=0$. Then it follows that $y\in\mathscr{I}_n$ (Lemma \ref{l:monsky}(c)). Since this holds for all $n$, Lemma \ref{l:inverselim} says $y=0$,
whence $x=0$ as desired.

Suppose $f\not=0$ and let $Z_1$,...,$Z_l$ be the $\Z_p$-flats in Theorem \ref{t:monsky}. For each $i$,
let $\chi(\xi_i)=\zeta_i$ be one of the defining equations \eqref{e:def} of $Z_i$, so that \eqref{e:Z} says $Z_i\subset \Delta_{\epsilon_i}$,
if $\zeta_i$ is of order $p^{r_i}$ and $\epsilon_i:=\nu_{\xi_i,r_i-1,r_{i}}$.  Then we put $\epsilon:=\epsilon_1\cdot\cdots\cdot \epsilon_l$
to have
\begin{equation}\label{e:epsilon}
\Delta_f\subset \Delta_\epsilon.
\end{equation}

Denote $\nu_{n,m}^{(i)}:=\nu_{\sigma_{i+1},n,m}\cdot\cdots\cdot \nu_{\sigma_d,n,m}$, $\underline{-1}=(-1,...,-1)$,
and
$\underline{m}^{(i)}=(m_1^{(i)},...,m_d^{(d)})$, with
$m_j^{(i)}=n$ for $j\leq i$; $m_j^{(i)}=m$, for $j>i$.  We claim that if
\begin{equation}\label{e:tibi}
\nu_{n,m}^{(i)} \cdot t_i\equiv b_i\cdot f \pmod{\mathscr{I}_{\underline{-1},\underline{m}^{(i)}}},
\end{equation}
holds for $t_i, b_i\in\Lambda_\Gamma$, then there is some $b_{i+1}\in\Lambda_\Gamma$ such that
$$
\nu_{n,m}^{(i+1)} \cdot \epsilon t_i\equiv b_{i+1}\cdot f \pmod{\mathscr{I}_{\underline{-1},\underline{m}^{(i+1)}}}.
$$
Then, beginning with \eqref{e:omegay}, for which $i=0$, by repeatedly applying the above implication, we can deduce
$$\epsilon^{d-1} y\equiv b_d\cdot f\pmod{\mathscr{I}_n},$$
which means $\omega\epsilon^{d-1} x\in \mathscr{I}_n+(f)$. As this happens for every $n$, Lemma \ref{l:inverselim} says
$\omega\epsilon^{d-1} x\in (f)$. The proposition is proved, if $f$ is relatively prime to $\omega\epsilon$, and
it remains to treat the case
where $f$ is a power of $\nu_{\sigma,\alpha-1,\alpha}$, $\alpha\geq 0$, $\sigma$ extendable to a $\Z_p$-basis of $\Gamma$,
because every irreducible factor of
$\omega\epsilon$ is of such type.
We may also assume that $\sigma=\sigma_1$.
Because $\tau_1,...,\tau_c$ form a tight generating
set, we can write
$$\omega=\nu_{\sigma,\alpha-1,\alpha}^\delta\cdot \omega',$$
with $\delta=0,1$, and $\omega'$ relatively prime to $f$. Note that if $\delta=1$, then $\alpha=0$ and $\sigma_1=\tau_j$, for some $j$, we may assume that the basis $\sigma_1,...,\sigma_d$
are taken from $\{\tau_1,...,\tau_c\}$. Put $z=\omega' x$, $\underline{\mathfrak a}=(\mathfrak a_1,...,\mathfrak a_d)$ with $\mathfrak a_1=\alpha$, $\mathfrak a_j=-1$ for $j\not=1$. By \eqref{e:x},
\begin{equation}\label{e:omega'z}
\nu_{n,m}\cdot z\equiv \omega' a\cdot f \pmod{\mathscr{I}_{\underline{\mathfrak a}, \underline{m}^{(0)}}}.
\end{equation}

Again, we claim that for $n>\alpha$, if
\begin{equation}\label{e:sici}
\nu_{n,m}^{(i)} \cdot s_i\equiv c_i\cdot f \pmod{\mathscr{I}_{\underline{\mathfrak a},\underline{m}^{(i)}}},
\end{equation}
holds for $s_i, c_i\in\Lambda_\Gamma$, then there is some $c_{i+1}\in\Lambda_\Gamma$ such that
$$
\nu_{n,m}^{(i+1)} \cdot s_i\equiv c_{i+1}\cdot f \pmod{\mathscr{I}_{\underline{\mathfrak a},\underline{m}^{(i+1)}}}.
$$
Then we can deduce that $\omega' x\in \mathscr{I}_{\underline{\mathfrak a}, \underline{n}^{(0)}}+(f)$, for all $n>\alpha$,
and hence
$\omega' x\in (f)$. Therefore, $\bar x=0$, since $\omega'$ is relatively prime to $f$.

The proofs of the claims rely on Lemma \ref{l:monsky}(c). Denote $\underline{\mathfrak b}^{(i)}:=(u_1,...,u_d)$, with $u_{i+1}=n$ and $u_j=-1$ for $j\not=i+1$.
Let $\chi\in \Delta_{\underline{\mathfrak b}^{(i)},\underline{m}^{(i)}}$. The inclusion
$$ \Delta_{\underline{\mathfrak b}^{(i)},\underline{m}^{(i)}}\subset \Delta_{\nu_{\sigma_{i+1},n,m}}\cap \Delta_{\underline{-1},\underline{m}^{(i)}}$$
leads to $\chi(\nu_{\sigma_{i+1},n,m})=0=\chi(\mathscr{I}_{\underline{-1},\underline{m}^{(i)}})$, so \eqref{e:tibi} yields $\chi(b_i\cdot f)=0$.
Now, if $\chi(f)=0$, then \eqref{e:epsilon} implies $\chi(\epsilon)=0$; otherwise, we have $\chi(b_i)=0$. Therefore, $\chi(\epsilon b_i)=0$
always holds. Lemma \ref{l:monsky}(c) says $\epsilon b_i\in \mathscr{I}_{\underline{\mathfrak b}^{(i)},\underline{m}^{(i)}}$.
Then we can write
$$\epsilon b_i=\nu_{\sigma_{i+1},n,m}\cdot b_{i+1}+b_{i+1}',$$
with $b_{i+1}'\in (\nu_{\sigma_1,-1,n},...,\nu_{\sigma_{i,-1,n}},\nu_{\sigma_{i+2},-1,m},...,\nu_{\sigma_d,-1,m})\subset \mathscr{I}_{\underline{-1},\underline{m}^{(i+1)}}$. By \eqref{e:tibi} again,
$$\nu_{\sigma_{i+1},n,m}\cdot (\nu_{n,m}^{(i+1)}\cdot \epsilon\cdot t_i-b_{i+1}\cdot f)\in \mathscr{I}_{\underline{-1},\underline{m}^{(i+1)}}.$$
In view of Lemma \ref{l:monsky}, this proves the first claim, because $\Delta_{\nu_{\sigma_{i+1},n,m}}\cap \Delta_{\underline{-1},\underline{m}^{(i+1)}}=\emptyset$,
if $\chi\in \Delta_{\underline{-1},\underline{m}^{(i+1)}}$, then $\chi( \nu_{\sigma_{i+1},n,m})\not=0$ and consequently
$\chi(\nu_{n,m}^{(i+1)}\cdot t_i-b_{i+1}\cdot f)=0$.

The proof of the second claims is similar to the previous one, the basic difference lies in
that all characters applied will be in $\Delta_{\nu_{\sigma, \alpha, \beta}}$, for some $\beta>\alpha$,
so that $\chi(f)$ is never zero. Let $\underline{\mathfrak a}^{(i)}
=(\mathfrak a_1^{(i)},
...,\mathfrak a_d^{(i)})$ be such that $\mathfrak a_j^{(i)}=\mathfrak a_j$, for $j\not=i+1$, $\mathfrak a_{i+1}^{(i)}=n$. Let
$\chi\in\Delta_{\underline{\mathfrak a}^{(i)},\underline m^{(i)}}$. Since $\chi(\nu_{\sigma_{i+1},n,m})=0
=\chi(\mathscr{I}_{\underline{\mathfrak a},\underline m^{(i)}})$, the congruence \eqref{e:sici} yields $\chi(c_i\cdot f)=0$,
hence $\chi(c_i)=0$. Lemma \ref{l:monsky}(c) says $ c_i\in \mathscr{I}_{\underline{\mathfrak a}^{(i)},\underline{m}^{(i)}}$,
and we can write
$$c_i=\nu_{\sigma_{i+1},n,m}\cdot c_{i+1}+c_{i+1}',$$
with $c_{i+1}'\in (\nu_{\sigma_1,\mathfrak a_1,n},...,\nu_{\sigma_{i,\mathfrak a_i,n}},\nu_{\sigma_{i+2},\mathfrak a_{d+2},m},...,\nu_{\sigma_d,\mathfrak a_d,m})\subset \mathscr{I}_{\underline{\mathfrak a},\underline{m}^{(i+1)}}$. Thus, by \eqref{e:sici},
$$\nu_{\sigma_{i+1},n,m}\cdot (\nu_{n,m}^{(i+1)}\cdot s_i-c_{i+1}\cdot f)\in \mathscr{I}_{\underline{\mathfrak a},\underline{m}^{(i+1)}}.$$
But $\chi(\nu_{\sigma_{i+1},n,m})\not=0$, for all $\chi\in\Delta_{\underline{\mathfrak a},\underline{m}^{(i+1)}}$,
it follows that $\nu_{n,m}^{(i+1)}\cdot s_i-c_{i+1}\cdot f$ is contained in $\mathscr{I}_{\underline{\mathfrak a},\underline{m}^{(i+1)}}$.
\end{proof}

\begin{remark}\label{r:rmk} If we replace $\mathscr{J}_n\fY$, $\mathscr{J}_m\fY$ in {\em{Definition \ref{d:ddot}}} by smaller
$\mathfrak{I}_n$,
$\mathfrak{I}_m$ satisfying $\nu_{n,m}(\mathfrak I_n)\subset \mathfrak I_m$, for $m>n$, and let $\Dot{\Dot{\fY}}_{(n)}$, $\Dot{\Dot{\fY}}$ be the resulting
counterpart of $\Ddot{\fY}_{(n)}$, $\Ddot{\fY}$, then $ \Dot{\Dot{\fY}}\subset \Ddot{\fY}\subset \fY$. Therefore, $ \Dot{\Dot{\fY}}\sim 0$ as well.
See the proof of {\em{Theorem \ref{t:x0}}}.
\end{remark}

\subsection{The module $\dot X_K$}\label{su:x0} 
Now we prove Theorem \ref{t:x0}.
By replacing $k$ by $k_n$ if necessary, we may assume that
for each $v$ in the ramification locus $S$, there is some integer $e$ such that
\begin{equation}\label{e:torfree}
\Gamma/\Gamma_v^0\simeq \Z_p^e,
\end{equation}
and
\begin{equation}\label{e:union}
\Gamma=\sum_{v\in S}\Gamma_v^0.
\end{equation}

Let $K_n$ be the maximal unramified abelian $p$-extension over $k_n$ so that by Class Field Theory,
$A_n:=A_{k_n}=\Gal(K_n/k_n)$. For $m>n$, the norm map $A_m\longrightarrow A_n$ is compatible with the
restriction of Galois action $\Gal(K_m/k_m)\longrightarrow \Gal(K_n/k_n)$, we have $X_K=\Gal(L/K)$, for
$L=\bigcup_n K_n$. Denote $G:=\Gal(L/k)$. Then $G/X_K=\Gamma$ and since $L/K$ is unramified,
at every place $v$ of $k$, we have $G_v^0\simeq \Gamma_v^0$, a commutative group.
$$
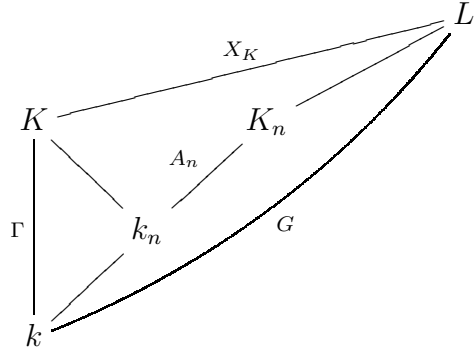
\xymatrix{
&  &   &  &{L} \\
{K}\ar@{-}[dd]_{\Gamma} \ar@{-}[urrrr]^{X_K} \ar@{-}[dr]& 
&{K_n} \ar@{-}[urr]& &\\
&{k_n}\ar@{-}[ur]^{A_n} & & &\\
{k}\ar@{-}[ur]  \ar@/_1.2pc/@{-}[uuurrrr]_{G}& & & &
}$$
Suppose $S=\{v_1,...,v_s\}$. For each $j$, choose a place $u_j$ of $L$ sitting
over $v_j$ and then choose a $\Z_p$-basis $\tilde \xi_1^{(j)},...,\tilde \xi_{d_j}^{(j)}$ of $G_{u_j}^0$. By \eqref{e:torfree} and \eqref{e:union},
we can choose these bases to have the union of their images under $G\longrightarrow\Gamma$ form a tight set
$\mathfrak g:=\{\tau_1,...,\tau_c\}$ of generators of $\Gamma$. Then among $\mathfrak g$, we choose a $\Z_p$-basis
$\{\sigma_1,...,\sigma_d\}$ of $\Gamma$. We lift each $\sigma_i$ to some $\tilde\xi_l^{(j)}$ and denote it by $\tilde\sigma_i$.
Every $g\in G$ can be uniquely written as
\begin{equation}\label{e:express}
g=\tilde\sigma_1^{a_1}\cdot\cdots\cdot \tilde\sigma_d^{a_d}\cdot x_g,\;\;a_1,...,a_d\in\Z_p, x_g\in X_K.
\end{equation}

Let $J_{ram}$ denote the $\Z_p$-submodule of $X_K$ generated by
$$\{x_g\;\mid\; g=\tilde\gamma\cdot \tilde\xi_i^{(j)}\cdot \tilde\gamma^{-1},\;\text{for some}\; j=1,...,s, i=1,...,d_j, \tilde\gamma\in G\}.$$
The commutator $\tilde\sigma_i\cdot\tilde\sigma_j\cdot\tilde\sigma_i^{-1}\cdot \tilde\sigma_j^{-1}\in J_{ram}$, because
$g=\tilde\sigma_i\cdot\tilde\sigma_j\cdot\tilde\sigma_i^{-1}$ in \eqref{e:express}
yields
$$\tilde\sigma_i\cdot\tilde\sigma_j\cdot\tilde\sigma_i^{-1}=\tilde\sigma_j\cdot x_{i,j},\;\; x_{ij}\in J_{ram}.$$

The canonical action of $\Gamma$ on $X_K$ coincides with the conjugation in $G$, namely, if $\gamma\in\Gamma$,
$x\in X_K$, and  $\tilde \gamma\in G$ is a lifting of $\gamma$, then
$$\tensor[^\gamma] x{}=\tilde\gamma \cdot x\cdot \tilde\gamma^{-1}.$$
Hence $\tensor[^{\gamma-1}] x{}$ is the same as the commutator $\tilde\gamma \cdot x\cdot \tilde\gamma^{-1}\cdot x^{-1}$ in $G$.
Set
$$J:=\mathscr{I}_0\cdot  X_K+J_{ram}\subset X_K.$$

Note that $J$ is a $\Gamma$-module, since $\Gamma$ acts trivially on $X_K/\mathscr{I}_0 X_K$.
Let $\tilde J\subset G$ denote the closed subgroup generated by $J$ and $\{\tilde\sigma_1,...,\tilde\sigma_d\}$.
Then for every place $u$ of $L$, the inertia subgroup $G_u^0\subset \tilde J$, because $u=\tensor[^{\gamma}]u{_j}$ for some $j$,
and hence $\tilde J$ contains the $\Z_p$-basis $\tilde\gamma\cdot \tilde\xi_i^{(j)}\cdot \tilde\gamma^{-1}$, $i=1,...,d_j$, of $G_u^0$.
Thus, $\tilde J$ is the closed subgroup of $G$ topologically generated by all commutators and all inertia subgroups of $G$.
Therefore, the fixed field $L^{\tilde J}$ is the maximal unramified abelian $p$-extension $K_0$ of $k$, with $\Gal(K_0/k)=A_0$.
By the uniqueness of the expression \eqref{e:express}, we have 
\begin{equation}\label{e:gj}
X_K/J\simeq G/\tilde J=\Gal(K_0/k)=A_0.
\end{equation}
Also, \eqref{e:union} implies $K\cap K_0=k$, hence $\Gal(KK_0/K)\simeq \Gal(K_0/k)$. By \eqref{e:gj},
\begin{equation}\label{e:KK0}
\Gal(KK_0/K)=X_K/J\simeq  A_0.
\end{equation}

To apply the above argument to $k_n$, let $G^{(n)}$ denote the pre-image of $\Gamma^{(n)}$ under $G\longrightarrow \Gamma$ and set
$\tilde\sigma_{i,n}:=\tilde\sigma_i^{p^n}$. Every element $g\in G^{(n)}$ can be uniquely expressed as
$$g=\tilde\sigma_{1,n}^{a_1}\cdot\cdots\cdot \tilde\sigma_{d,n}^{a_d}\cdot x_g,\;\;a_1,...,a_d\in\Z_p, x_g\in X_K.$$

Write $\tilde\xi_{1,n}^{(j)},...,\tilde\xi_{d_j,n}^{(j)}$ for $(\tilde\xi_1^{(j)})^{p^n},...,(\tilde\xi_{d_j}^{(j)})^{p^n}$.  They form a $\Z_p$-basis of
the inertia subgroup of $G_{u_j}^{(n)}$.
Let $J_{ram}^{(n)}$ denote the $\Z_p$-submodule of $X_K$ generated by
$$\{x_g\;\mid\; g=\tilde\gamma\cdot \tilde\xi_{i,n}^{(j)}\cdot \tilde\gamma^{-1},\;\text{for some}\; j=1,...,s, i=1,...,d_j, \tilde\gamma\in G\}.$$
Set
$$J_n:=\mathscr{I}_n\cdot  X_K+J_{ram}^{(n)}\subset  X_K.$$
Let $\tilde J_n\subset G$ denote the closed subgroup generated by $J_n$ and $\{\tilde\sigma_{1,n},...,\tilde\sigma_{d,n}\}$.

\begin{lemma}\label{l:yn} The fixed field of $\tilde J_n\subset G^{(n)}$ is $K_n$.
We have
$$\Gal(KK_n/K)=X_K/J_n\simeq G^{(n)}/\tilde J_n=\Gal(K_n/k_n)=A_n.
$$
Furthermore, $J_n$ is a $\Gamma$-module, hence a $\Lambda_\Gamma$-module.
\end{lemma}
\begin{proof}
It remains to show that $J_n$ is invariant under the action of $\Gamma$, or equivalently $K_n/k$ is a Galois extension.
But this follows from the fact that $k_n/k$ is Galois and $K_n$ is the maximal unramified abelian $p$-extension
over $k_n$.
\end{proof}

The next lemma  follows \cite[Theorem 7]{iwa73}.
\begin{lemma}\label{l:nu}
If $m\geq n$, then $\nu_{n,m}(J_n)\subset J_m$ and that induces the commutative diagram
$$\xymatrix{X_K/J_n \ar[r]^-\sim \ar[d]^-{\nu_{n,m}} & \;A_n\ar[d]^-{c_{k_m/k_n}}\\
X_K/J_m \ar[r]^-\sim  & \,\;A_m.}$$
\end{lemma}
\begin{proof} 
Since \eqref{e:union} says $K/k$ contains no non-trivial unramified subextension, the restriction of Galois action $X_K\longrightarrow A_n$ is surjective for every $n$.  Let $x_n\in A_n$ and let
$x\in X_K$ such that $x\mid_{k_n}=x_n$.  Denote $x\mid_{k_m}:=x_m$. Let $[\;]_m$ denote the Artin map and let $\mathfrak L$ be
an ideal in $\O_{k_m}$ with $[\mathfrak L]_m=x_m$.  If $\mathfrak l=\Nm_{k_m/k_n}(\mathfrak L)$, then $x_n=[\mathfrak l]_n$, hence
$$c_{k_m/k_n}(x\mid_{k_n})=c_{k_m/k_n}(x_n)=[ \Nm_{k_m/k_n}(\mathfrak L)]_m=\prod_{\gamma\in\Gal(k_m/k_n)} \tensor[^\gamma] {[\mathfrak L]}{_m}=\nu_{n,m}(x)\mid_{k_m}.$$
In particular, if $x\in J_n$, then the left-hand side is trivial, hence $\nu_{n,m}(x)\in J_m$.
\end{proof}

\begin{proof}[Proof of {\em{Theorem \ref{t:x0}}}] By Lemma \ref{l:nu}, we can write for the capitulation kernel
$$\dot {X}_n=\bigcup_n \ker (\xymatrix{ X_K/J_n \ar[r]^{\nu_{n,m}} &  X_K/J_n}).$$
This makes it possible to apply the technique developed in the previous sections. Because $J\subset X_K$ is of finite index,
for convenience, we replace $X_K$ by $J$. To be more precise, since $\mathscr{I}_n\cdot X_K\subset J_n$, we have $\mathscr{I}_n\cdot J\subset J_n$, put $\dot J_{(n,m)}:=\ker (\,\xymatrix{J/J_n \ar[r]^-{\nu_{n,m}} & J/J_m}\,)$,
$\dot J_{(n)}:=\bigcup_{m\geq n} \dot J_{(n,m)}$, and $\dot J:=\varprojlim_n \dot J_{(n)}$.
Then it is sufficient to show $\dot J\sim0$.

Take $\fY:=J$, $\mathfrak I_n:=J_n\cap \mathscr{J}_n\cdot \fY$ and let $\Dot{\Dot\fY}_{(n)}$, $\Dot{\Dot\fY}$ be as in Remark \ref{r:rmk}. The homomorphisms $\fY/\mathfrak I_n\longrightarrow \fY/\mathscr{J}_n\cdot\fY$, for all $n$,
yield a homomorphism $\Dot{\Dot\fY}\longrightarrow \Ddot \fY$ fitting into the commutative diagram
$$\xymatrix{\Dot{\Dot\fY} \ar[r]\ar[d] & \Ddot \fY\ar[d]\\
        \fY \ar@{=}[r] & J.}$$
Since the down-arrows are injective, Proposition \ref{p:monsky} implies $\Dot{\Dot\fY}\sim 0$.

To see the difference between $\mathfrak I_n$ and $J_n$, we first check that for $x\in X_K$,
the element $\omega_{\sigma_i,n}\cdot x=\nu_{\sigma_i,0,n}\cdot \omega_{\sigma_i,0}\cdot x\in \mathscr{J}_n\cdot J$, because $\sigma_i\in \mathfrak g$ and
$\omega_{\sigma_i,0}\cdot x\in \mathscr{I}_0\cdot X_K\subset J$. This shows $\mathscr{I}_n\cdot X_K\subset
\mathscr{J}_n\cdot J\cap J_n=\mathfrak I_n$. Denote $\tilde\xi_i^{(j)}=:g$, $\tilde\xi_{i,n}^{(j)}=:g_n$ and let $\tau\in\mathfrak g$
be the image of $g$ under $G\longrightarrow \Gamma$. Let $x=x_g\in J_{ram}$ so that $g=\bar{g}\cdot x$ with $\bar{g}=\tilde\sigma_1^{a_1}\cdot\cdots\cdot\tilde\sigma_d^{a_d}$. Then
$$\tilde\sigma_{1,n}^{a_1}\cdot\cdots\cdot\tilde\sigma_{d,n}^{a_d}\cdot x_{g_n}=g_n=g^{p^n}=(\bar g\cdot x)^{p^n}=\bar{g}^{p^n}\cdot \nu_{\tau,0, n} x.$$
Similarly, if $g'=\tilde \gamma\cdot g \cdot\tilde\gamma^{-1}=\bar g\cdot x'$, with $x'=x_{g'}\in J_{ram}$, and
$g'_n=\tilde \gamma\cdot g_n \cdot\tilde\gamma^{-1}$, then
$$\tilde\sigma_{1,n}^{a_1}\cdot\cdots\cdot\tilde\sigma_{d,n}^{a_d}\cdot x_{g_n'}=g_n'=(g')^{p^n}=(\bar g\cdot x')^{p^n}=\bar{g}^{p^n}\cdot \nu_{\tau,0, n}
x'.$$
Let $\tilde\rho:=(\tilde\sigma_{1,n}^{a_1}\cdot\cdots\cdot\tilde\sigma_{d,n}^{a_d})^{-1}\bar{g}^{p^n}$ and write $X_K$ additively. Then,
since $x_{g_n}, x_{g_n'}\in J_{ram}^{(n)}$,
$$x_{g_n}-x_{g_n'}=\tensor[^{\tilde\rho}]{(\nu_{\tau,0,n}(x-x'))}{}\in \mathfrak I_n.$$
Therefore, we have shown that the $\Z_p$-module $J_n/\mathfrak I_n$ is generated by the classes of $x_{g_n}$, with $g=\tilde\xi_i^{(j)}$,
$j=1,...,s$, $i=1,...d_j$. Then we observe that $x_{g_n}$ is trivial if the above $\tau=\sigma_i$. This shows the $p$-rank of
$J_n/\mathfrak I_n$ is at most $c-d$, which equals $0$ if $d=1$. The exact sequence
$$0\longrightarrow J_n/\mathfrak I_n\longrightarrow \fY/\mathfrak I_n\longrightarrow J/J_n\longrightarrow 0$$
gives rise to the exact sequence
$$\mathfrak a_n\longrightarrow \Dot{\Dot\fY}_{(n)}\longrightarrow \dot J_{(n)} \longrightarrow \mathfrak b_n,$$
with $p$-ranks of both $\mathfrak a_n$ and $\mathfrak b_n$ bounded by $c-d$. Consequently, the cokernel 
of the induced map
$\Dot{\Dot\fY}\longrightarrow \dot J$ is of $p$-ranks bounded by $c-d$. Hence $\dot J$ is pseudo-null and the theorem is proved for
$\dot X_K$. 
The corresponding assertion for $\dot X_K'$ can be proved by similar approach, of which we only give a sketch. 
Write 
$$\P=\P_1\sqcup \P_2,$$
where $\P_1$ consists of place $v$ such that 
$\Gamma_w^0\not=\Gamma_w$. 
For each $v\in \P$, choose a place
$u$ of $L$ above $v$ and denote $w=u\mid_{_K}$, if $v\in S$, let $u$ be as before.
If $v\in\P_1$, then by \eqref{e:torfree}, we have non-canonically 
$\Gamma_w=Z_p\times \Gamma_w^0$, whence $\Gamma_w\simeq G_u$ and $w$ splits completely over $L$. 
In this case, choose an $\eta_u\in G_u$ representing a topological generator of $G_u/G_u^0$ such that the map $G\longrightarrow \Gamma$
sends $\eta_u$ into $\mathfrak g$. Write $\eta_{u,n}$ for
$\eta_u^{p^n}$ and let $J_{un}^{(n)}$ denote the
$Z_p$-module generated by
$$\{x_g\;\mid\; g=\tilde\gamma\cdot \eta_{u,n}\cdot \tilde\gamma^{-1},\;\; v\in\P_1, \;\tilde\gamma\in G\}.$$
If $v\in \P_2$, then we have canonically $G_u\simeq\Gal(L_u/K_w)\times \Gamma_w$ with $\Gal(L_u/K_w)\subset X_K$, $G_u^0\simeq\Gamma_w$. Let $[X_K]$ be the $\Gamma$-submodule of $X_K$ generated by $\Gal(L_u/K_w)$, for all $v\in\P_2$.  
Set 
$$\{J_n\}=J_n+J_{un}^{(n)}+[X_K],$$
let $\widetilde{\{J_n\}}$ denote the closed subgroup of $G$ generated by $\{J_n\}$ and $\{\tilde\sigma_{1,n},...,\tilde\sigma_{d,n}\}$,
and let $K_n'$ be the fixed field of $\widetilde{\{J_n\}}$. Then $K_n'$ is the maximal unramified $p$-extension of $k_n$ with all places
above $p$ splitting completely, so $\Gal(K_n'/k_n) =A_n'$. Let $[A]_n$ be the image of $\{J_n\}$ under $X_K\longrightarrow A_n=X_K/J_n$.
Then the exact sequences 
$$0\longrightarrow [A]_n\longrightarrow A_n \longrightarrow A_n'\longrightarrow 0$$ 
and the fact that
$$[X_K]=\bigcap_n \{J_n\}=\varprojlim_n [A]_n$$ 
yield $X_K'=X_K/[X_K]$. Let $J_n'$ be the image of $\{J_n\}$ under $X_K\longrightarrow X_K'$. 
Then $(X'_K,J_n',K_n')$-version of Lemma \ref{l:yn} holds, namely,
$$X'_K/J'_n\simeq \Gal(K_n'/k_n)=A_n'.$$
Again by the calss field theory, the $(X'_K,J'_n,A'_n)$-version of Lemma \ref{l:nu} holds.
Then take $\fY':=J_0'$, $\mathfrak I_n':=J_n'\cap \mathscr{J}_n\cdot \fY'$ and proceed as above.
\end{proof}


\section{Cohomology groups of global units}\label{s:gu}
The proof of Theorem \ref{t:m} involves cohomology of unit groups.
Denote the group of global units of $E$ by $U_E:=\O_E^*$ and put $U_K:=\bigcup_E U_E$.
The cohomology groups $\coh^i(\Gamma, U_K)$, $i=1,2$, has been studied in
\cite{iwa83,yam84} for $d=1$ case. We are going to show that for general $d$, they are co-finitely generated $\Z_p$-modules.
Put
$$\mathcal U_E:=\Q_p/\Z_p\otimes U_E,\;\; \mathcal U_K:=\varinjlim_E \mathcal U_E$$ 
and
$$\mathfrak M_E:=\Gal(\bar{E}^{ab,p}/E),\;\; \mathfrak M_K:=\varprojlim_E \mathfrak M_E,$$
where $\bar{E}^{ab,p}$ is the maximal pro-$p$ abelian extension of $E$, unramified outside $p$.

\subsection{An exact sequence}\label{su:exact}
Let $\Div_E$ and $P_E$ denote the groups of divisors and principal divisors of $E$.
In \cite[Proposition 1]{iwa83}, Iwasawa deduces the seven-term exact sequence:
\begin{equation}\label{e:seven}
\xymatrix{0\ar[r] & \coh^1(G, U_E) \ar[r] & \Div_E^G/P_k \ar[r]^-{\alpha_{E/k}} & C_E^{G} \ar[r]^-{\beta_{E/k}} & \coh^2(G,U_E)
\ar[dll]_-{\gamma_{E/k}}\\
 &  & b_E  \ar[r] & \coh^1(G,  C_E) \ar[r] & \coh^3(G, U_E)  ,}
\end{equation}
with $G:=\Gal(E/k)$ and
$ b_E:=\ker (\coh^2(G,E^*)\longrightarrow \coh^2(G, \Div_E))$.

The exact sequence
identifies $\coh^1(G,U_E)$ with the kernel of $\alpha_{E/k}$. The restriction of $\alpha_{E/k}$ to $ C_k\subset \Div_E^G/P_k$ is
the capitulation homomorphism $c_{E/k}$. Therefore, we have the exact sequence
\begin{equation}\label{e:coh1}
\xymatrix{0\ar[r] & \ker (c_{E/k}) \ar[r] & \coh^1(G,U_E) \ar[r] & (\Div_E^G\cap P_E)/(\Div_k\cap P_E)\ar[r] & 0.}
\end{equation}

If $v\in S$ and $w$ is a place of $E$ sitting over $v$, then the sum
$$\epsilon_{v} :=\sum_{\gamma\in G/G_v} \tensor[^\gamma] w{}$$
is fixed by $G$. Set $g_{E/k,v}:=|G_v^0|$. In $\Div_E$, we have $g_{E/k,v}\cdot \epsilon_v=v$. The sequence
\begin{equation}\label{e:divdiv}
\xymatrix{0\ar[r] & \Div_k \ar[r] & \Div_E^G \ar[r]^-{l_{E/k}} & \bigoplus_{v\in S}  (\Z_p/g_{E/k,v}\Z_p) \cdot \epsilon_v\ar[r] & 0, }
\end{equation}
with $l_{E/k}$ the projection onto the $S_E$-component, is exact. Combining \eqref{e:coh1} and \eqref{e:divdiv}, we obtain the exact sequence
\begin{equation}\label{e:coh1e}
\xymatrix{0\ar[r] & \ker (c_{E/k}) \ar[r] & \coh^1(G,U_E) \ar[r]^-{\mathfrak l_{E/k}} & \bigoplus_{v\in S}  \Z_p/g_{E/k,v}\Z_p .} 
\end{equation}
By taking the injective limit of the above sequence, one can show that $\coh^1(\Gamma, U_K)$ is co-finitely generated over $\Z_p$.
As for $\coh^2(\Gamma,U_K)$, we have the following. Denote 
$$A_K:=\varinjlim_E A_E.$$

\begin{lemma}\label{l:equiv}
The following statements are equivalent:
\begin{enumerate}
\item[(a)] The $\Lambda_\Gamma$-module $W_K$ is finitely generated.
\item[(b)] The abelian group $A_K^\Gamma$ has finite $p$-rank.
\item[(c)] The abelian group $\coh^2(\Gamma, U_K)$ has finite $p$-rank.
\end{enumerate}
\end{lemma}
\begin{proof}
Since $A_K^\Gamma[p]$ is the Pontryagin dual of $W_K/(p\Lambda_\Gamma+\mathscr I_\Gamma)W_K$,
the equivalence between (a) and (b) follows from Nakayama's Lemma.

Denote $Q_E:=\image (\alpha_{E/k})\cap  A_E^\Gamma$. Then \eqref{e:seven} induces the exact sequence
$$\xymatrix{Q_K \ar[r] & A_K^\Gamma \ar[r] & \coh^2(\Gamma, U_K) \ar[r]^-\beta & \coh^2(k, K^*)}.
$$
Now, the $p$-rank of $Q_K$ is $\leq$ the $p$-rank of $A_k$ plus the cardinality of $S$. Also,
if $v\not\in S$, then the composition of $\xymatrix{\coh^2(\Gamma, U_K) \ar[r]^-\beta & \coh^2(k, K^*)\ar[r]
& \coh^2(k_v,K_v^*)}$ is the trivial homomorphism. Hence $\image(\beta)$ is embedded into
$\prod_{v\in S} \Q_p/\Z_p$ by the local invariant maps and is of $p$-rank bounded by the cardinality of $S$. Therefore, (b) and (c) are
equivalent. \end{proof}

\subsection{The structure of $\mathcal U_K^\vee$}\label{su:ulvee} 
Recall that $U_K$ denotes the group of global units of $K$ and 
$\mathcal U_K^\vee=\Hom(\Q_p/\Z_p\otimes U_K, \Q_p/\Z_p)$.

Let $\mathcal V_E$ denote the group of $\P_E$-units of $E$.
Kummer's theory yields the commutative diagram
$$\xymatrix{
U_E\otimes \frac{1}{p^n}\Z_p/\Z_p \ar[r] \ar@{^(->}[d] & \V_E \otimes \frac{1}{p^n}\Z_p/\Z_p \ar[r] \ar@{^(->}[d] &
E^*\otimes \frac{1}{p^n}\Z_p/\Z_p \ar[d]^-\simeq \\
\coh^1(\mathfrak M_E,\Bmu_{p^n}) \ar@{=}[r]&  \coh^1(\mathfrak M_E,\Bmu_{p^n}) \ar@{^(->}[r]  & \coh^1(k,\Bmu_{p^n}).
   }
$$
Since $E^*/\V_E$ and $\V_E/U_E$ are torsion free, arrows in the upper row of the diagram are injective.
Let $n$ goes to $\infty$. We obtain the commutative diagram
\begin{equation}\label{e:kummer}
\xymatrix{\mathcal U_E \ar@{^{(}->}[r] \ar[d]^-{i_E} & \V_E\otimes\Q_p/\Z_p \ar@{^{(}->}[r] \ar[d]^-{j_E} &  \coh^1(\mathfrak M_E,\Bmu_{p^\infty})\ar[d]^-{r_E}\\
\mathcal U_K  \ar@{^{(}->}[r] & \V_K\otimes \Q_p/\Z_p\ar@{^{(}->}[r] &  \coh^1(\mathfrak M_K,\Bmu_{p^\infty}).}
\end{equation}

\begin{lemma}\label{l:control}
The restriction map $\xymatrix{\coh^1(\mathfrak M_E,\Bmu_{p^\infty})\ar[r]^-{r_E} &
\coh^1(\mathfrak M_K,\Bmu_{p^\infty})^{\Gamma_E}}$ has finite kernel and cokernel.
\end{lemma}
\begin{proof}
\cite[Lemma 3.2.1]{tan10}.
\end{proof}

\begin{lemma}\label{l:controlv}
The cokernel of $\xymatrix{\V_E\otimes\Q_p/\Z_p \ar@{^{(}->}[r]  &  \coh^1(\mathfrak M_E,\Bmu_{p^\infty})}$
is finite.
\end{lemma}
\begin{proof}
We show that as
$n\rightarrow\infty$,
the corkernel of
$$\xymatrix{\V_E \otimes \frac{1}{p^n}\Z_p/\Z_p\ar@{^{(}->}[r] & \coh^1(\mathfrak M_E,\Bmu_{p^n})}$$
remains bounded. By Kummer's theory, an element of $\coh^1(\mathfrak M_E,\Bmu_{p^n})$ is represented by an $f\in E^*$ modulo
$(E^*)^{p^n}$ such that $p^n\mid \ord_v(f)$ for every $v\notin \P_E$.
Let $w_1,...,w_r$ be a set of generators of $A_E$ and let $p^m$ be the exponent of $A_E$
so that each $w_i^{p^m}=(f_i)$ for some $f_i\in E^*$. If $n>m$, then modulo $ (E^*)^{p^n}$,
$f$ can be expressed as a product of powers of $f_1^{p^{n-m}},...,f_r^{p^{n-m}}$ together with elements of $\V_E$. Hence the
cokernel in question has order bounded by $p^{mr}$.
\end{proof}

\begin{corollary}\label{c:control}
The natural map $\mathcal U_E\longrightarrow \mathcal U_K^{\Gamma_E}$ has finite kernel, its cokernel is cofinitely generated over $\Z_p$ of corank bounded by $|S_E|$.
\end{corollary}
\begin{proof} The exact sequence
$$\xymatrix{U_E \ar@{^{(}->}[r] & \V_E \ar[r]^-{\tilde\mho_E} & \prod_{v\in \P_E} \Z,}$$
where $\tilde\mho_E$ is defined by taking valuations at all $v\in \P_E$ induces the commutative diagram of complexes
$$\xymatrix{\mathcal U_E \ar@{^{(}->}[r] \ar[d]^-{\alpha_E} & \V_E\otimes\Q_p/\Z_p \ar[r]^-{\mho_E} \ar[d]^-{\beta_E}  &
\prod_{v\in \P_E} \Q_p/\Z_p \cdot v\ar[d]\\
\mathcal U_K  \ar@{^{(}->}[r] & \V_K\otimes \Q_p/\Z_p\ar[r]^-{\mho_K} &  \prod_{w\in \P_K} \Q_p/\Z_p \cdot w.}
$$
Write $Q_E:=\ker (\mho_E)$, $R_E:=\ker (\mho_K\circ \beta_E)$. Then $Q_E/\mathcal U_E$ is finite.
In view of the diagram \eqref{e:kummer} and Lemma \ref{l:control}, \ref{l:controlv}, we need to show that the cokernel of
$\xymatrix{Q_E \ar@{^{(}->}[r] & R_E}$ has corank bounded by $|S_E|$. If $w$ is a place of $K$ sitting over $v$ of $E$ with
$v\not\in S_E$, then $w$ is unramified under $K/k$, and hence the map
$$\Q_p/\Z_p\cdot v\longrightarrow \Q_p/\Z_p\cdot w$$
is injective. Therefore, we have the exact sequence
$$\xymatrix{0\ar[r] & Q_E \ar[r] & R_E \ar[r] &\prod_{v\in S_E} \Q_p/\Z_p\cdot v.}$$
\end{proof}

\begin{lemma}\label{l:rankq} Let $\mathfrak Q_K$ be a finitely generated $\Lambda_\Gamma$-module.
If $\mathfrak Q_K$ is of rank $r$ over $\Lambda_\Gamma$, then $\mathfrak Q_K/\mathscr I_n \mathfrak Q_K$ has $\Z_p$-rank
$$r_n=r\cdot p^{dn}+O(p^{(d-1)n}).$$
\end{lemma}
\begin{proof} There exists an exact sequence
$$\xymatrix{0\ar[r] & \bigoplus_{i=1}^r \Lambda_\Gamma \ar[r] & \mathfrak Q_K \ar[r] & T \ar[r] & 0.}$$
Since $T$ is torsion, \cite[Lemma 1.13]{tan14} says $T/\mathscr I_nT$ has $\Z_p$-rank $=O(p^{(d-1)n})$. Hence
$$r_n\leq r\cdot p^{dn}+O(p^{(d-1)n}).$$
To obtain an inequality in the other direction, we use an exact sequence
$$\xymatrix{  \mathfrak Q_K \ar[r] & \bigoplus_{i=1}^r \Lambda_\Gamma \ar[r] & T' \ar[r] & 0,}$$
with $T'$ torsion.
\end{proof}

Let $r_1$ and $r_2$ denote the number of
real and complex places of $K$. 

\begin{proposition}\label{p:calu}
The $\Lambda_\Gamma$-module $\mathcal U_K^\vee$ is finitely generated, of rank $r_1+r_2$.
\end{proposition}

\begin{proof}
Since $|S_{k_n}|=O(p^{(d-1)n})$, Corollary \ref{c:control} says
$\mathcal U_K^\vee/\mathscr I_n \mathcal U_K^\vee$, which is the Pontryagin dual of  $\mathcal U_K^{\Gamma^{(n)}}$,
has $\Z_p$-rank equals $(r_1+r_2)\cdot p^{dn}+O(p^{(d-1)n})$. In particular, $\mathcal U_K^\vee/(p\Lambda_\Gamma+\mathscr I_\Gamma)\mathcal U_K^\vee$ is finite. Hence $\mathcal U_K^\vee$ is finitely generated. By Lemma \ref{l:rankq}, the rank of $\mathcal U_K^\vee$
equals $r_1+r_2$. 
\end{proof}

\subsection{The module $W_K$}\label{su:ml}

Let $K'/k$ be an intermediate $\Z_p^{d-1}$-extension of $K/k$ with Galois group $\Gamma'$.
We shall fix a $\Z_p$-extension $k^{(\infty)}/k$ linearly disjoint from $K'/k$ so that $K=K'k^{(\infty)}$ and $K'\cap k^{(\infty)}=k$.
Let $k^{(n)}$ denote the $n$th layer of $k^{(\infty)}/k$. From now on, the symbol $F$ will denote a finite intermediate extension of $K'/k$, while $E$ denotes a finite intermediate extension of $K/k$ such that
$E\cap K'=F$ and $E=Fk^{(n)}$ for some $n$. Put $F^{(\infty)}:=Fk^{(\infty)}$.
Denote $\Psi=\Gal(K/K')$ and let $\psi$ be a topological generator of $\Psi$.

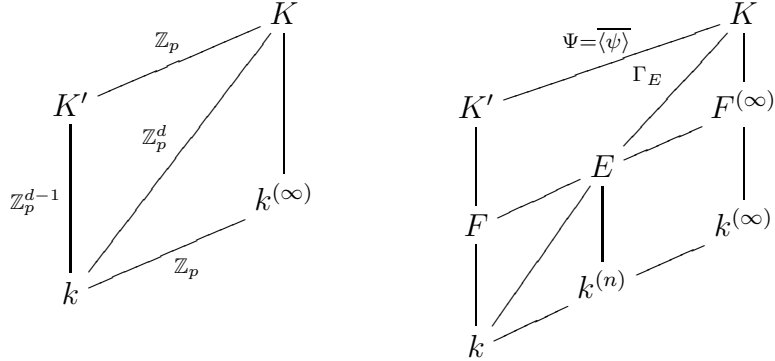
$$\xymatrix@R=.5pc{ & &K\\
& & \\
{K'}  \ar@{-}[rruu]^{\mathbb{Z}_p} & &\\
& & \\
 & &{k^{(\infty)}} \ar@{-}[uuuu]\\& & \\
 k \ar@{-}[uuuu]^{\mathbb{Z}_p^{d-1}}
 \ar@{-}[rruu]_{\mathbb{Z}_p}
 \ar@{-}[rruuuuuu]^{\mathbb{Z}_p^{d}} & &
}
\qquad \qquad
\xymatrix@R=.5pc{ & &K\\
& & \\
{K'}  \ar@{-}[rruu]^-{\Psi=\overline{\langle \psi\rangle}} & &{F^{(\infty)}}  \ar@{-}[uu]\\
 &{E} \ar@{-}[ruuu]^{\Gamma_E} \ar@{-}[ru] &\\
 F \ar@{-}[ru] \ar@{-}[uu] & &{k^{(\infty)}} \ar@{-}[uu] \\
   &{k^{(n)}} \ar@{-}[uu] \ar@{-}[ru]&\\
   k \ar@{-}[uu] \ar@{-}[ru] \ar@{-}[ruuu] & &
}$$

For each $F$ and for $m\geq n$, we have the commutative diagram
$$\xymatrix{(\Z[\frac{1}{p^n}]/\Z)\otimes U_F \ar[r]^-{p^n}_-\sim \ar[d] &  U_F/p^n U_F \ar@{->>}[r] \ar[d]^-{p^{m-n}} 
&\coh^2(Fk^{(n)}/F, U_{Fk^{(n)}})\ar[d]^-{inf}\\ 
(\Z[\frac{1}{p^m}]/\Z)\otimes U_F \ar[r]^-{p^m}_-\sim &  U_F/p^m U_F \ar[r] & \coh^2(Fk^{(m)}/F, U_{Fk^{(m)}}),}
$$
where the left down-arrow is induced by the natural map $\Z[\frac{1}{p^n}]\longrightarrow \Z[\frac{1}{p^m}]$ and the right down-arrow
is the inflation map. This yields 

\begin{equation}\label{e:nk}
\xymatrix{\mathcal U_F \ar@{->>}[r] & \coh^2((F^{(\infty)}/F, U_{F^{(\infty)}}).}
\end{equation}

\vskip5pt
\begin{proof}[Proof of Theorem {\em{\ref{t:m}}}] Since $W_K'\subset W_K$, it is sufficient to treat $W_K$.
Lemma \ref{l:equiv} says
the theorem is equivalent to the assertion that $\coh^2(\Gamma, U_K)$ has finite $p$-rank.
For $d=1$, the value of this $p$-rank is shown in
\cite{iwa83, yam84}, the equality \eqref{e:nk} also implies it
is bounded by the rank of $U_k$. We prove the theorem by induction on $d$.

By \eqref{e:seven}, we have the exact sequence
$$\xymatrix{0\ar[r]& \dot{C}_{E/F} \ar[r] &  C_E^\Psi \ar[r] & \coh^2(E/F,U_E),} $$
where $\dot{C}_{E/F}$ is generated by $\Div_E^\Psi$. Denote $\mathcal C_K:=\varinjlim C_E$. Let $E$ goes to $K$ and take the direct limit to obtain
the exact sequence
$$\xymatrix{0\ar[r]& \dot{\mathcal C}_{K/K'} \ar[r] & \mathcal C_K^\Psi \ar[r] & \coh^2(\Psi,U_K).} $$
Since the $p$-part of $\mathcal C_K^\Gamma$ is Pontryagin dual to $W_K/\mathscr I_\Gamma W_K$,
it remains to show both $\dot{\mathcal C}_{K/K'}^{\Gamma'}$ and $\coh^2(\Psi,U_K)^{\Gamma'}$ are of finite $p$-ranks.

The surjection  \eqref{e:nk} says $\coh^2(\Psi,U_K)$ is a quotient of $\mathcal U_{K'}$. Hence its Pontryagin dual, denoted $\mathfrak J_{K'}$,
is a $\Lambda_{\Gamma'}$ submodule
of $\mathcal U_{K'}^\vee$.  Proposition \ref{p:calu} implies $\mathfrak J_{K'}$ is $\Lambda_{\Gamma'}$ finitely generated.
Then $\coh^2(\Psi,U_K)^{\Gamma'}$, being the Pontryagin dual of $ \mathfrak J_{K'}/\mathscr I_{\Gamma'}\mathfrak J_{K'}$,
must have finite $p$-rank.

Let $T\subset S$ be the subset consisting of places $v$ with $\Gamma_v^0$ of rank $d$. We choose $K'$ such that $K/K'$ is unramified outside
$T_{K'}$, which is a finite set. Apply \eqref{e:divdiv} to the case where $k=F$ and then let $E$ goes to $K$ to obtain the exact sequence
$$\xymatrix{0\ar[r] & \bar{\mathcal C}_{K'} \ar[r] & \dot{\mathcal C}_{K/K'} \ar[r] & V\ar[r] & 0,}$$
where $\bar{\mathcal C}_{K'}$ is the image of the capitulation homomorphism $\mathcal C_{K'}\longrightarrow \mathcal C_K$ and $V$ is a quotient of $\prod_{w\in T_{K'}}\Q_p/\Z_p\cdot w$. Obviously, the $p$-rank of $V$ is finite. Now $\bar{\mathcal C}_{K'}^\vee$ is a $\Lambda_{\Gamma'}$ submodule of $W_{K'}$, whence finitely generated by the induction hypothesis. Consequently,
$\bar{\mathcal C}_{K'}^{\Gamma'}$, the Pontryagin dual of $\bar{\mathcal C}_{K'}^\vee/\mathscr I_{\Gamma'}\bar{\mathcal C}_{K'}^\vee$,
must have finite $p$-rank. The above exact sequence implies the $p$-rank of $\dot{\mathcal C}_{K/K'}^{\Gamma'}$ is finite.
\end{proof}

Write $\mathfrak a_E:=A_E$, $\mathfrak b_E:=\Hom(A_E,\Q_p/\Z_p)$. For $E\subset E'\subset K$, let $\mathfrak r^{E'}_E:\mathfrak b_E\longrightarrow \mathfrak b_{E'}$
and $\mathfrak k^{E'}_E:\mathfrak b_{E'}\longrightarrow \mathfrak b_{E}$ denote respectively
the restriction and the corestriction.  Extend them to homomorphisms
$$\mathfrak r^{E'}_E:\mathfrak a_E\times \mathfrak b_E\longrightarrow \mathfrak a_{E'}\times \mathfrak b_{E'}$$
and
$$\mathfrak k^{E'}_E:\mathfrak a_{E'}\times \mathfrak b_{E'}\longrightarrow \mathfrak a_{E}\times \mathfrak b_{E}$$
such that the restrictions to the first factors are respectively the capitulation and the norm. Then
$$\mathfrak r^{E'}_E\circ \mathfrak k^{E'}_E=\Nm_{\Gal(E'/E)},$$
the norm map of $\Gal(E'/E)$-action.
Let $\langle \;,\;\rangle_E:\mathfrak a_E\times \mathfrak b_E\longrightarrow \Q_p/\Z_p$ be the natural pairing.
For $(a,b) \in \mathfrak a_E\times\mathfrak b_E$, $(a',b') \in \mathfrak a_{E'}\times\mathfrak b_{E'}$, we have
$$\langle a, \mathfrak k^{E'}_E(b')\rangle_E=\langle \mathfrak r^{E'}_E(a),b'\rangle_{E'},\;\;\;
\langle a', \mathfrak r^{E'}_E(b)\rangle_{E'}=\langle \mathfrak k^{E'}_E(a'),b\rangle_{E}.$$
Thus, in the terminology of \cite[\S 3.1.1]{lltt18}, the collection
$$\mathfrak A:=\{\mathfrak a_E, \mathfrak b_E,\langle \;,\;\rangle_E,\mathfrak r^{E'}_E,\mathfrak k^{E'}_E\;\mid\; k\subset E \subset E'\subset K\}$$
form a complete $\Gamma$-system. Indeed, it is a
{\bf{T}}-system [op. cit.], because for every intermediate $\Z_p^e$-extension $M/k$ of $K/k$, $d\geq e\geq 1$, if
$$\mathfrak a_M:=\varprojlim_{k\subset E\subset M} \mathfrak a_E,\;\;\; \mathfrak b_M:=\varprojlim_{k\subset E\subset M} \mathfrak b_E,$$
then $\mathfrak a_M\times \mathfrak b_M$ is finitely generated torsion over the Iwasawa algebra of $M/k$.


\vskip5pt
\begin{proof}[Proof of Theorem {\em{\ref{t:mw}}}] In view of  \cite[Theorem 1]{lltt18}, for proving $W_K\sim X_K^\sharp$,
we need to show the system $\mathfrak A$ is pseudo-control, in the sense that
$$\dot{\mathfrak a}_K\times \dot{\mathfrak b}_K:=\varprojlim_E \bigcup_{E\subset E'} \ker(\mathfrak r^{E'}_E)$$
is pseudo-null. That $\dot{\mathfrak a}_K\sim 0$ is by Theorem \ref{t:x0}. To show $\dot{\mathfrak b}_K\sim 0$,
we may assume that $K/k$ has no nontrivial unramified extension. In this case $\xymatrix{\mathfrak b_E\ar[r]^-{\mathfrak r^{E'}_E} & \mathfrak b_{E'}}$ is injective, whence $\dot{\mathfrak b}_K=0$.

The assertion $W'_K\sim X_K'^\sharp$ is proved similarly. We first form the $\Gamma$-system $\mathfrak A'$ with $\mathfrak a_E'=A_E'$,
$\mathfrak b_E'=\Hom(\mathfrak a_E',\Q_p/\Z_P)$. Since $\mathfrak a_M'$ is quotient of $\mathfrak a_M$, $\mathfrak b_M'\subset \mathfrak b_M$,
$\mathfrak A'$ is a $T$-system.  Theorem \ref{t:x0} says $\dot{\mathfrak a}_K'\sim 0$, and since
$\dot{\mathfrak b}_K'\subset \dot{\mathfrak b}_K\sim 0$, $\mathfrak A'$ is pseudo-controlled. Then apply  \cite[Theorem 1]{lltt18}.
\end{proof}

As for \eqref{e:wchi},  form the $\Gamma$-system $\mathfrak A_\chi=\{\mathfrak a_{\chi,E},\mathfrak b_{\chi,E},\mathfrak r_E^{E'},\mathfrak k_E^{E'}\}$ with $\mathfrak a_{\chi,E}=(A_E)_{\chi^{-1}}$,
$\mathfrak b_{E,\chi}=\Hom(\mathfrak a_{\chi,E},\Q_p/\Z_P)$, it is a pseudo-controlled $T$-system, and hence the first pseudo-isomorphism follows.
The second is proved similarly.

\end{document}